\documentclass{article}
\usepackage{latexsym}
\usepackage{amssymb}
\usepackage{amsmath}%
\newtheorem{theorem}{Theorem}
\newtheorem{lemma}[theorem]{Lemma}
\newtheorem{corollary}[theorem]{Corollary}
\newcommand\QED{\ifhmode\allowbreak\else\nobreak\fi\quad\nobreak$\Box$\medbreak}
\newenvironment{proof}{\par\noindent{\bf Proof:}}{\rm\enspace{\QED\par}}
\newcommand\commentout[1]{}

\begin{document}
\begin{titlepage}
\title{The Expected Order of a Random Unitary Matrix \\
(Preliminary Version)}
\author{
Eric Schmutz\\
\small Department of Mathematics\\
\small Drexel University\\
\small Philadelphia, Pa. 19104 \\
\small Eric.Jonathan.Schmutz@drexel.edu
 }

\date{\today}
\maketitle
\abstract{
Let $U(n,q)$ be the group consisting of those invertible matrices $A=\left(a_{i,j}\right)_{1\leq i,j\leq n}$ whose inverse is the conjugate transpose with respect to the involution
$c\mapsto c^{q}$ of the finite field ${\mathbb F}_{q^{2}}$. In other words, the $i,j$'th entry
of $A^{-1}$ is $a_{j,i}^{q}.$
Let $\mu_{n}=\frac{1}{|U(n,q)|}\sum\limits_{A\in U(n,q)}{\rm Order}(A)$ be the average of the orders of the elements in this finite group. We prove the following conjecture of Fulman:  for any fixed $q$, as $n\rightarrow\infty$,
$$\log\mu_{n}={n}\log (q) -\log n+o_{q}(\log n).$$
}
\vskip.5cm \noindent
{\bf Keywords and phrases}: {\em  Unitary group, finite field, cycle index }
\end{titlepage}
\section{Introduction}
\par
This paper concerns the finite unitary group $U(n,q)$, so we begin by reviewing some basic notation and definitions
relating to this group.
Let $q=p^{\ell}$ for some prime number $p$ and some positive integer $\ell$.
  The involution $c\mapsto c^{q}$ is an automorphism of the finite field ${\mathbb F}_{q^{2}}$ that fixes
the subfield ${\mathbb F}_{q}$. If  $A=\left(a_{i,j}\right)_{1\leq i,j\leq n}$ is an $n\times n$ matrix with entries in 
${\mathbb F}_{q^{2}}$, let $A^{*}$ be the matrix whose $i,j$'th entry is  $a_{j,i}^{q}$ (for $1\leq i,j\leq n).$
Define the  unitary group $U(n,q)$ to be the group  consisting of those $n\times n$ matrices $A$
 for which $A^{-1}=A^{*}.$   
It is 
well known (e.g.\cite{Grove}, page 109) that, under  matrix multiplication, this set of matrices forms a group of order
\begin{equation}
\label{orderunq}
|U(n,q)|=q^{n^{2}}\prod\limits_{j=1}^{n}(1-\frac{(-1)^{j}}{ q^{j}}).
\end{equation}

\par
For any prime power $r$, let $GL(n,r)$ be the group of invertible $n\times n$ matrices with entries in ${\mathbb F}_{r}$.
It is well known  that $GL(n,r)$ has order 
\begin{equation}
\label{ordergln}
|GL(n,r)|=r^{n^{2}}\prod\limits_{j=1}^{n}(1-\frac{1}{ r^{j}}).
\end{equation}
Note that $U(n,q)$ is a subgroup of $GL(n,q^{2})$, not $GL(n,q)$.

\par

For any finite group $G$, let $\mu(G)=\frac{1}{ |G|}\sum\nolimits_{g\in G}{\bf V}(g)$, where ${\bf V}(g)$ is the order of $g$.  
Stong \cite{St1}proved that, for any prime power $r$, 
\begin{equation}
\label{gln}
\log \mu(GL(n,r))=n\log r-\log n+o_{r}(\log n).
\end{equation}
as $n\rightarrow\infty.$
 Fulman proposed  the analogous 
problem estimating 
$\mu(U(n,q))$. He proved that
$\log \mu(U(n,q))\geq \frac{1}{ 2}{n}\log q^{2} -\log n+o_{q}(\log n),$  
and conjectured that this lower bound is sharp insofar as the 
\lq\lq $\geq$\rq\rq can be replaced with \lq\lq$=$\rq\rq.
The goal of this paper is to prove Fulman's conjecture. 
\par
 The rest of this section contains additional  definitions and symbols that are listed  in quasi-alphabetical order, and then used globally
 without comment. 
\begin{itemize}
\item $|f|$: the degree of the polynomial $f$.
\item $\tilde{f}$: if $f(x)=x^{d}+\sum\limits_{j=0}^{d-1}a_{j}x^{j}$
 is a monic polynomial of degree $d$ with non-zero constant term $a_{0}$,
then $\tilde{f}(x)=x^{d}+ \sum\limits_{j=0}^{d-1}a_{0}^{-q}a_{d-j}^{q}x^{j}.$
\item $\lbrack\hskip-1pt\lbrack z^{n}\rbrack\hskip-1pt\rbrack  F(z)$: coefficient of $z^{n}$ in $F(z).$
\item $c_{i}(\pi)$: number ofparts of size $i$ that the partition $\pi$ has.
\item $C_{\infty}=\prod\limits_{j=1}^{\infty}(1-\frac{1}{ 2^{j}})\approx.289$
\item $E_{n}$: expected value with respect to $P_{n}$, i.e. for any real valued function  ${\bf Y}$ that is defined on 
characteristic polynomials of matrices in $U(n,q)$,
        $E_{n}({\bf Y})=\frac{1}{ |U(n,q)|}
\sum\limits_{A\in U(n,q)}{\bf Y}(char.poly.(A)).$
\item
${\cal I}_{d,r}$ =set of all monic polynomials of degree $d$ in $\mathbb{ F}_{r}[x]$ that are irreducible over ${\mathbb F}_{r}$ (except for $\phi(x)=x,$ which is excluded from ${\cal I}_1.$). 
\item ${\cal I}_{d}={\cal I}_{d,q^{2}}$
\item
${\cal I}=\bigcup\limits_{d=1}^{\infty}{\cal I}_{d,q^{2}}.$
\item  ${\cal J}_{d}$ =monic,  irreducible polynomials $\phi$ of degree $d$ in $\mathbb{ F}_{q^{2}}[x]$  that satisfy $\phi=\tilde{\phi}$ 
\item ${\cal J}=\bigcup\limits_{d=1}^{\infty}{\cal J}_{d}.$
\item ${\cal K}_{d}={\cal I}_{d,q^{2}}-{\cal J}_{d}=$
monic,  irreducible polynomials $\phi$ of degree $d$ in $\mathbb{ F}_{q^{2}}[x]$  that satisfy $\phi\not=\tilde{\phi}$ 
\item ${\cal K}=\bigcup\limits_{n=1}^{\infty}{\cal K}_{d}$
\item ${\cal K}_{+},{\cal K}_{-}:$  disjoint subsets of ${\cal K}$ such that $\phi \in {\cal K}_{+}$iff $\tilde{\phi} \in {\cal K}_{-}.$
\item
$m_{\phi}=m_{\phi}(f)=$ the multiplicity of $\phi$ in $f$:
for $\phi\in {\cal I}$ and $f\in \mathbb{F}_{q^{2}}[x]$, 
$\phi^{m_{\phi}(f)}$ divides $f$ but $\phi^{m_{\phi}(f)+1}$ does not divide $f$.
\item $m_{\phi}(A)=m_{\phi}($characteristic polynomial of  $A).$
\item ${\bf M}=\max\limits_{\phi\in {\cal I}}m_{\phi}.$
\item $Q_{b}$:set of all partitions of $b$ into distinct odd parts.
\item $\Omega_{n}$= all characteristic polynomials of matrices in $U(n,q)$=  monic, degree $n$,
 polynomials $f\in {\mathbb F}_{q^{2}}[x]$ satisfying $m_{\phi}(f)=m_{\tilde{\phi}}(f)$ for all $\phi\in {\cal I}.$

\item ${\bf P}_{n}=$ the probability measure on $\Omega_{n}$ that is induced by the uniform distribution on $U(n,q),$
i.e. ${\bf P}_{n}({\cal S})=\frac{|\bigl\lbrace A\in U(n,q): char.poly(A)\in {\cal S}\bigr\rbrace|}{ |U(n,q)|}$ for all ${\cal S}\subseteq \Omega_{n}.$
\item $Q_{n}=$ set of all partitions of $n$ into distinct parts.
\item $\tau_{\phi}$=order of the roots of the irreducible polynomial $\phi$ (as multiplicative units in
in the splitting field for $\phi$.)
\item ${\bf T}(f)=LCM\lbrace \tau_{\phi}:\phi $ is an in irreducible factor of $f\rbrace$
\item ${\bf X}_{1}(f)=LCM(\bigl\lbrace q^{|\phi|}+1: \ m_{\phi}(f)>0, \phi\in {\cal J}\bigr\rbrace).$
\item ${\bf X}_{2}(f)=LCM(\bigl\lbrace q^{2|\phi|}-1: \ m_{\phi}(f)>0, \phi\in {\cal K}_{+}\bigr\rbrace).$
\item ${\bf X}_{1}(\pi)=LCM(\bigl\lbrace q^{d}+1: \pi$ has a part of size $d \bigr\rbrace).$
\item ${\bf X}_{2}(\lambda)=LCM(\bigl\lbrace q^{2d}-1: \lambda$ has a part of size $d \bigr\rbrace).$
\item ${\bf X}(f)=$ least common multiple of  ${\bf X}_{1}(f)$ and  ${\bf X}_{2}(f).$
\item  ${\bf X}(A)={\bf X}($ characteristic polynomial of $A).$
\item ${\bf V}(A)$ =order of $A = \min\bigl\lbrace e: A^{e}=I\bigr\rbrace \rbrace.$
\end{itemize}
\par

\section{Reduction from ${\bf V} $ to ${\bf X}$}
There is a close relationship between the order of a matrix $A\in GL(n,q^{2})$ and the 
orders of its eigenvalues (as multiplicative units in a splitting field for the 
characterstic polynomial). Hence we begin this section with a   simple lemma about the orders of 
the roots of irreducible polynomials. We also state, for future reference, Fulman's  formula for the number of unitary matrices with a given characteristic polynomial.
These facts are used to bound the maximum order, and
to prove that most matrices in $U(n,q)$ do not have eigenvalues of large algebraic multiplicity.  This in turn enables us to reduce the problem of estimating
$E_{n}({\bf V})$ to the easier problem of estimating $E_{n}({\bf X}).$ 
\par
Recall that, if 
 $\phi(x)=x^{d}+\sum\limits_{j=0}^{d-1}a_{j}x^{j}$
 is a monic polynomial of degree $d$ with non-zero constant term $a_{0}$, then
 $\tilde{\phi}(x)=x^{d}+ \sum\limits_{j=0}^{d-1}a_{0}^{-q}a_{d-j}^{q}x^{j}.$ 
\begin{lemma}\label{Mullen} Suppose $\phi\in {\cal I}_{d},$ 
and suppose $\tau_{\phi}$ and $\tau_{\tilde{\phi}}$ are respectively the orders of the roots of $\phi$  and $\tilde{\phi}$ (as multiplicative units in ${\mathbb F}_{q^{2d}}).$  Then
\begin{itemize}
\item   $\tau_{\phi}=\tau_{\tilde{\phi}}$ 
\item If $\phi=\tilde{\phi}$, then  $\tau_{\phi}$ is a divisor of $ q^{d}+1.$
\end{itemize}
\end{lemma}

\begin{proof}
Observe that
$\rho$ is a root of $\phi$ if and only if $\rho^{-q}$ is a root of  $\tilde{\phi}$:
\begin{eqnarray}
\tilde{\phi}(\rho^{-q})=a_{0}^{-q}\rho^{-dq}\sum\limits_{k=0}^{d}a_{k}^{q}\rho^{kq}&\\
=a_{0}^{-q}\rho^{-dq}\left(\sum\limits_{k=0}^{d}a_{k}\rho^{k}\right)^{q}.
\end{eqnarray}
As an element of ${\mathbb F}_{{q}^{2d}}^{*},$  the order of $\rho^{-q}$ 
is equal to the order of 
its inverse  $\rho^{q}$ , which is in turn equal to  the order of $\rho$ (since $q$ and 
$q^{2d}-1$ are coprime.).  This proves the first part:  $\tau_{\phi}=\tau_{\tilde{\phi}}.$ 
\par
Let $\rho$ be one of the roots of $\phi$, assume that
$\phi=\tilde{\phi}.$   
Then $\rho^{-q}$ must be one of the roots of $\phi$. 
But the roots of $\phi$
are $\rho^{q^{2}},\rho^{q^{4}},\dots ,\rho^{q^{2d-2}}, \rho^{q^{2d}}=\rho.$
Hence, for some positive integer $j\leq d$, we have $\rho^{q^{2j}}=\rho^{-q}$, and consequently
$\rho^{q(q^{2j-1}+1)}=1.$  This proves that $\tau_{\phi}$ divides $q(q^{2j-1}+1)$.
But $\tau_{\phi}$ also divides $q^{2d}-1$,  and $g.c.d.(q,q^{2d}-1)=1$. 
Therefore
$\tau_\phi$  divides $q^{2j-1}+1.$ 
Let $m$ be the smallest positive integer such that $\tau_{\phi}$ divides $q^{m}+1.$
If $\tau_{\phi}=2$, then it is clear that $\tau_{\phi}$ divides $q^{d}+1$ since $\tau_{\phi}$
divides $q^{2d}-1=(q^{d}+1)(q^{d}-1)$ and both factors are even. We may therefore assume that $\tau_{\phi}>2$.
Using   Proposition 1 of  \cite{YucasMullen} (with $s=2d$),  and the fact that $\tau_{\phi}| q^{2d}-1$, we get
  $2d=2\ell m$ for some positive integer $\ell$.
  We know $d$ is odd (Fulman \cite{Fulman1}, Theorem 9),
therefore $\ell$ must also be odd.  Again using  Proposition 1 of  \cite{YucasMullen} (this time with $s=d$), we get $\tau_{\phi}|q^{d}+1.$
\end{proof}
Let $\Omega_{n}$ be the set of polynomials that are characteristic polynomials of matrices in
$U(n,q).$ A beautiful characterization of these polynomials is known. 
A monic polynomial $f$ is in $\Omega_{n}$ if and only if $m_{\phi}(f)=m_{\tilde{\phi}}(f)$ for all $\phi\in {\cal I}$; the multiplicity of $\phi$ is the same as the multiplicty of $\tilde{\phi}$
for all irreducible polynomials $\phi$.
In fact, with  the notational convention that $|U(0,r)|=|GL(0,r)|=1$ for all 
prime powers $r$, we can state the following theorem of
 Fulman\cite{Fulman1}:
\begin{theorem}
\label{charpoly}
(Fulman) If $f\in \Omega_{n},$  then
$$P_{n}(\lbrace f\rbrace )=\prod\limits_{\phi\in {\cal J}}
\frac{
q^{|\phi|(m_{\phi}^{2}-m_{\phi})} }{ |U(m_{\phi},q^{|\phi|})|
}\cdot
\prod\limits_{\theta\in {\cal K}_{+}}
\frac{
q^{2|\theta|(m_{\theta}^{2}-m_{\theta})}}{ |GL(m_{\theta},q^{2|\theta|})|}
$$
\end{theorem}
Theorem \ref{charpoly} was just one application of powerful generating function techniques
that Fulman developed for $U(n,q)$ and other finite classical groups.
Related work can be found in Kung\cite{Kung}, Stong\cite{St2}, and recent work of 
Fulman, Neumann, and  Praeger, e.g. \cite{FPN}.

\par
  If the eigenvalues are all distinct, then 
the order of a matrix is just the least common multiple of the orders of the eigenvalues. 
The general case is a bit more complicated because the Jordan form includes off-diagonal elements. 
This leads to Theorem \ref{VMTineq} below. This convenient  inequality 
is an immediate consequence of the slightly stronger inequality in 
the introduction of Stong's paper \cite{St1}.
(See also Lidl and Niederreiter\cite{Lidl}, page 80):
\begin{theorem} 
\label{VMTineq} 
For all  $ A\in GL(n,{q^{2}}),$ 
$ {\bf V}\leq  p{\bf M} {\bf T}.$
\end{theorem}
An immediate consequence of Theorem \ref{VMTineq} is a bound on the maximum order:
\begin{corollary}
For all $A\in GL(n,q^{2}),$ ${\bf V}< pnq^{2n}.$
\end{corollary}
However a stronger inequality holds for $U(n,q).$
\begin{corollary}
\label{maxv}
For all $A\in U(n,q),$ ${\bf V} \leq  3p{\bf M}q^{n}.$
\end{corollary}

\begin{proof}
By Theorem \ref{VMTineq}, 
it suffices to prove that ${\bf T} \leq 3q^{n}.$
Suppose the characteristic polynomial of $A$ is
$$\prod\limits_{i=1}^{r}\phi_{i}^{m_{\phi_{i}}}
\prod\limits_{j=1}^{s}(\phi_{r+j}\tilde{\phi}_{r+j})^{m_{\phi_{r+j}}}$$
where $\phi_{i}\in {\cal J}$ for $i\leq r$ and $\phi_{r+j}\in {\cal K}_{+}$ for $j\leq s.$
To simplify notation, let $d_{i}=|\phi_{i}|$, and $\tau_{i}=\tau_{\phi_{i}}$.
 Then by Lemma \ref{Mullen}, $\tau_{i}$ divides $(q^{d_{i}}+1)$ for $i\leq r$ and 
 $\tau_{i}$ divides $q^{2d_{i}}-1$ for $r< i\leq r+s$. Hence
\begin{align}
&{\bf T}(A)=LCM( \tau_{1},\tau_{2}\dots \tau_{r+s})\leq\\
\label{repeats}
&\leq \cdot LCM(q^{d_{1}}+1,q^{d_{2}}+1)\dots ,q^{d_{r}}+1)LCM( q^{2d_{r+1}}-1,\dots ,q^{2d_{r+s}}-1)
\end{align}
Without loss of generality, assume $d_{i}\not=d_{j}$  for $1\leq i<j\leq r.$
(If two degrees are equal, then we can  remove one of the arguments to the least common multiple function without changing its value.)
Then
\begin{align}
&{\bf T}(A)
\leq \prod\limits_{i=1}^{r}(q^{d_{i}}+1)\cdot \prod\limits_{j=1}^{s}q^{2d_{r+j}}\\
&=q^{n}\prod\limits_{i=1}^{r}(1+\frac{1}{ q^{d_{i}}})\\
&\leq q^{n}\prod\limits_{i=1}^{r}(1+\frac{1}{ 2^{{i}}})\leq 3q^{n}.
\end{align}
\end{proof}

We have a bound on the maximum order, but we still need to prove that the maximum multiplicity ${\bf M}$
 is usually small. If
 $\xi=\xi(n)\rightarrow\infty$, then
with high probability,
no irreducible factor has multiplicity larger than $\xi$.  

\begin{lemma}
\label{multiplicitybound} For all positive integers $n,$ and all $\xi >2,$ 
$P_{n }({\bf M}>\xi)\leq 40q^{1-\xi}.$
\end{lemma}
\begin{proof} 
Suppose $d$ is a positive integer $\leq n$ and  $\psi=\tilde{\psi} \in J_{d}$.   Note that, for $f\in \Omega_{n}$,
we have $m_{\psi}(f)=\ell$ if and only if $f=\psi^{\ell}g$ for some $g\in \Omega_{n-d\ell}$
such that $m_{\psi}(g)=0$.
Hence, by Theorem \ref{charpoly},
\begin{eqnarray}
\label{pmeqell}
P_{n}(m_{\psi}=\ell)=\frac{q^{d\ell^{2}-d\ell}}{ |U(\ell,q^{d})|}P_{n-d\ell}(m_{\psi}=0) \leq \frac{q^{d\ell^{2}-d\ell}}{|U(\ell,q^{d})|} .
&\\
\end{eqnarray}
 Using (\ref{orderunq}), we get 
\begin{eqnarray}
\frac{q^{d\ell^{2}}}{ |U(\ell,q^{d})|}=\frac{1}{ \prod\limits_{j=1}^{\ell}
(1-\frac{(-1)^{j}}{ q^{dj}})}&\\
< \frac{1}{ \prod\limits_{j=1}^{\ell}(1-\frac{1}{ q^{dj}})} <\frac{1}{ C_{\infty}}<4.&\\
\end{eqnarray}
 Putting this back into the right side of  (\ref{pmeqell}), and summing on $\ell$, we get
\begin{equation}
\label{J}
P_{n}(m_{\psi}\geq\xi)=\sum\limits_{\ell\geq \xi}P_{n}(m_{\psi}=\ell)\leq 4\sum\limits_{\ell\geq\xi}q^{-\ell d} \leq 8q^{-d\xi}.
\end{equation}
Similarly, for any $\psi\in {\cal K}_{d},$ we have
\begin{eqnarray}
P_{n}(m_{\psi}=\ell)=\frac{q^{2d(\ell^{2}-\ell)}}{ |GL(\ell,q^{2d})|}P_{n-2d\ell}(m_{\psi}=0)&\\
\leq \frac{q^{2d(\ell^{2}-\ell)}}{ |GL(\ell,q^{2d})|}=
\frac{q^{-2d\ell}}{ \prod\limits_{j=1}^{\ell}(1-\frac{1}{ q^{2dj}})}&\\
< \frac{q^{-2d\ell}}{ \prod\limits_{j=1}^{\infty}(1-\frac{1}{ 2^{2j}})} < 2q^{-2d\ell},\\
\end{eqnarray}
and consequently
\begin{equation}
\label{notJ}
P_{n}(m_{\psi}\geq\xi) \leq 4q^{-2d\xi}.
\end{equation}
Now, given a real number $\xi>2$,   let ${\bf N}_{\xi}$ be the number of irreducible factors having multiplicity greater than $\xi$. 
Then ${\bf M}> \xi $ if and only if ${\bf N}_{\xi}>0,$ and it suffices to show that 
$P_{n}({\bf N}_{\xi}>0)\leq 40q^{1-\xi}.$ 

Combining (\ref{J}) and (\ref{notJ}),   we get
\begin{eqnarray}
P_{n}({\bf N}_{\xi}>0)\leq E({\bf N}_{\xi})=\sum\limits_{d=1}^{\lfloor n/d\xi\rfloor }\sum\limits_{\phi\in {\cal I}_{d,q^{2}}}P_{n}(m_{\phi}>\xi)&\\
\label{forjdbound}
=\sum\limits_{d=1}^{\lfloor n/d\xi\rfloor }\left(\sum\limits_{\phi\in {\cal J}_{d}}P_{n}(m_{\phi}>\xi)+ \sum\limits_{\phi\in {\cal K}_{d}}P_{n}(m_{\phi}>\xi)\right)          &\\
\label{rhs}
\leq \sum\limits_{d=1}^{\infty}\left(|{\cal J}_{d}|8q^{-d\xi}+|{\cal K}_{d}|4q^{-2d\xi}\right).
\end{eqnarray}
It is well known (e.g. \cite{Berl}, page 80) that, for any prime power $r$,
\begin{equation}
\label{Mignotte}
|{\cal I}_{d,r}|=\frac{1}{d}\sum\limits_{k|d}\mu(k)r^{d/k}\leq \frac{r^{d}}{ d}
\end{equation}
Since ${\cal K}_{d}\subseteq {\cal I}_{d,q^{2}}$, we follows that
\begin{equation}
\label{Kd}
|{\cal K}_{d}|\leq  \frac{q^{2d}}{d}.\end{equation}
We need a similar estimate for $|{\cal J}_{d}|.$ 
Fulman proved that 
\begin{equation}
\label{mobius}
|{\cal J}_{d}|=
\begin{cases}0& \text{if $d$ is even,}\\
\frac{1}{ d}\sum\limits_{k|d}\mu(k)(q^{d/k}+1)&\text{ else.}
\end{cases}
\end{equation}
It is well known that, for all $d>1$, $\sum\limits_{k|d}\mu(k)=0.$
Therefore, for all odd $d>1$, 
\begin{equation}
\label{Jd}
|{\cal J}_{d}| =\frac{1}{d}\sum\limits_{k|d}\mu(k)q^{d/k}=|{\cal I}_{d,q}|\leq \frac{q^{d}}{ d}.
\end{equation}
(It is interesting that  $|{\cal J}_{d}|$ is {\sl exactly equal } to 
$|{\cal I}_{d,q}|$ , even though the two sets are not equal.)
 For $d=1$ we have $|{\cal J}_{d}|=q+1\leq 2q$, so for all $d\geq 1$ we crudely have 
\begin{equation}
\label{Jd1}
|{\cal J}_{d}|\leq  \frac{2q^{d}}{ d}.
\end{equation}
For $0<x< \frac{1}{ 2},$ we have $-\log(1-x)<2x,$  and for $\xi>2,$ we have $q^{1-\xi}< \frac{1}{ 2}.$  Therefore, by
putting (\ref{Jd1}) and (\ref{Kd}) into (\ref{rhs}), we get
\begin{eqnarray}
P_{n}({\bf N}_{\xi}>0)\leq\sum\limits_{d=1}^{\infty}\left(\frac{16q^{d-d\xi}}{ d}+\frac{4q^{2d-2d\xi}}{ d}\right)&\\
\leq -20\log(1-q^{1-\xi})\leq 40q^{1-\xi}.
\end{eqnarray}

\end{proof}

Now that Lemma \ref{multiplicitybound} is available, we can
reduce the problem from the task of estimating  $E_{n}({\bf V})$ to the slightly easier task of 
estimating $E_{n}({\bf X}).$

\begin{lemma}
$\log E_{n}({\bf V})\leq \log E_{n}({\bf  X})+O(\log\log n).$
\label{finalreduction}
\end{lemma}
\begin{proof}
By lemma \ref{Mullen}, ${\bf T}$(A) divides ${\bf X}(A)$ for all $A$. It therefore suffices to prove that
\begin{equation}
\log E_{n}({\bf V})\leq \log E_{n}({\bf  T})+O_{p}(\log\log n).
\end{equation}

For any $\xi$, we have
\begin{equation}
\label{Treduction}
E_{n}({\bf V})= P_{n}({\bf M}\leq \xi )E_{n}({\bf V} | {\bf M}\leq \xi)+P_{n}({\bf M} >\xi)E_{n}({\bf V}|{\bf M}>\xi).
\end{equation}
To estimate the second term of the two terms on the right side of in (\ref{Treduction}), we use Corollary \ref{maxv} and Lemma \ref{multiplicitybound} 
with $\xi=\log^{2} n$:
\begin{equation}
\label{term2}
P_{n}({\bf M}>\xi)E_{n}({\bf V}|{\bf M}>\xi)\leq (40q^{1-\log^{2}n})(3pnq^{n})=
 q^{n-\log^{2}n(1+o(1))}.
\end{equation}
For the first term on the right side of  (\ref{Treduction}), we are conditioning on ${\bf M}\leq \xi$ so  
we can use the inequality
${\bf M\leq \xi}$ together with the inequality
${\bf V}\leq  p{\bf M} {\bf T}$  from  Theorem \ref{VMTineq}:
\begin{align}
&P_{n}({\bf M}\leq \xi )E_{n}({\bf V} | {\bf M}\leq \xi)\leq \\
p\xi & P_{n}({\bf M}\leq \xi )
 E_{n}({\bf T} | {\bf M}\leq \xi)\\
\leq p\xi &\left( P_{n}({\bf M}\leq \xi )
E_{n}({\bf  T} | {\bf M}\leq \xi)+P_{n}({\bf M}> \xi )E_{n}({\bf  T} | {\bf M}> \xi)\right)\\
=p\xi &E_{n}({\bf T}).
\label{term1}
\end{align}

Finally, putting (\ref{term1})and (\ref{term2}) back  into (\ref{Treduction}), we get
\begin{equation}
\label{prelog}
E_{n}({\bf V})\leq (p\log^{2} n)E_{n}({\bf  T})\left(1+ \frac{q^{n-\log^{2}n(1+o(1))}}{E_{n}({\bf T})}\right).
\end{equation}
Since $E_{n}({\bf T})\geq q^{n-\log n}$ for all sufficiently large $n$ (section 6 of Fulman \cite{Fulman1}),
the lemma follows from (\ref{prelog}) by taking logarithms.
\end{proof}

\section{ Key Factorization.}

There is a second factorization of characteristic polynomials that is crucial for this paper. The idea is to factor 
the characteristic polynomial  $f$ as $f=gh$ where
\begin{itemize}
\item {\bf X}(f)={\bf X}(g)
\item $g$ is easier to work with than $f$, and
\item $g$ and $h$ are themselves characterstic polynomials of unitary matrices.
\end{itemize}

To that end, define ${\cal D}(f)$ to be the set of polynomials $g$ that satisfy the following three conditions:
\begin{itemize}
\item $g\in \Omega$
\item $g$ divides $f$
\item ${\bf X}(g)={\bf X}(f)$
\end{itemize}
The set   ${\cal D}(f)$ is non-empty since $f\in {\cal D}(f).$   Because ${\cal D}(f)$ is a non-empty finite set that is  partially ordered by divisibility, we  can choose a minimal element $\pi(f).$  

Suppose we have chosen, for each  $f\in \Omega_{n}$,  a factor $g=\pi(f)$ that is minimal in  ${\cal D}(f).$
It is clear that, no matter how the minimal element is chosen, it will have the following useful properties:
\begin{itemize}
\label{observations}
\item For all $\phi$ in ${\cal I}$, $m_{\phi}(g)=0$ or $1$.
\item For all positive integers $d$, $\pi(f)$ has zero,one, or two  irreducible factors of degree $d$. If there is one such irreducible factor $\phi$, then $\phi\in {\cal J}_{d}.$ If here are two, and $\phi$ is one of them, then $\tilde{\phi}$ is the other and both are in $ {\cal K}_{d}.$
\end{itemize}

The third  property we need is less obvious, but it is proved in the following lemma.
\begin{lemma}
\label{sub}
If $f\in \Omega_{n}$ and $g=\pi(f)$ has degree $|g|<n$, and if $h=\frac{f}{ \pi(f)}$, then
$$P_{n}(\lbrace f\rbrace )\leq P_{|g|}(\lbrace g\rbrace )P_{n-|g|}(\lbrace h\rbrace ).$$
\end{lemma}
\begin{proof}
We consider each factor of $P_{n}(\lbrace f\rbrace )$ in the factorization of Theorem \ref{charpoly}
show that it is bounded  above by the corresponding factors in the product  $P_{|g|}(\lbrace g\rbrace )P_{|h|}(\lbrace h\rbrace ).$

\par
Suppose first that 
 $\phi\in{\cal J}_{d}$ for some $d$, and suppose $\phi$ divides $g$. To simplify notation, let $m=m_{\phi}(f).$
In Theorem \ref{charpoly}, the factor of $P_{n}(\lbrace f\rbrace )$  corresponding to $\phi$  is
\begin{eqnarray}
\frac{q^{d(m^{2}-m)}}{ |U(m,q^{d})|}=\frac{q^{-dm}}{ \prod\limits_{j=1}^{m}(1-\frac{(-1)^{j}}{ q^{dj}})} &\\
=\frac{q^{-d}}{ (1-\frac{(-1)^{m}}{ q^{dm}})}\frac{q^{-d(m-1)}}{ \prod\limits_{j=1}^{m-1}(1-\frac{(-1)^{j}}{ q^{dj}})} &\\
\label{AB}
\leq \frac{q^{-d}}{(1-\frac{1}{ q^{d}})}\frac{q^{-d(m-1)}}{ \prod\limits_{j=1}^{m-1}(1-\frac{(-1)^{j}}{ q^{dj}})} 
\end{eqnarray}
Since $\phi$ divides $g$, we have  $m_{\phi}(g)=1$ and $m_{\phi}(h)=m-1.$
 Therefore the  factor of $P_{|g|}(g)$ that corresponds to $\phi$ is 
	$\frac{q^{-d}}{(1-\frac{1}{ q^{d}})}$, and the factor of $P_{|h|}(h)$ that corresponds to $\phi$ is
$\frac{q^{-d(m-1)}}{ \prod\limits_{j=1}^{m-1}(1-\frac{(-1)^{j}}{ q^{dj}})}.$ These are precisely the two factors on the right of (\ref{AB}).
\par
Similarly, if $\phi\in {\cal K}^{+}$ has degree $d$ and $\phi$ divides $g$,
then the factor of $P_{n}(\lbrace f\rbrace )$ that corresponds to $\phi$ is
\begin{eqnarray}
\frac{q^{2d(m^{2}-m)}}{ |GL(m,q^{2d})|}=\frac{q^{-2dm}}{ \prod\limits_{j=1}^{m}(1-\frac{1}{ q^{2dj}})} &\\
\label{AB2}
\leq \frac{q^{-2d}}{ (1-\frac{1}{ q^{2d}})}\frac{q^{-2d(m-1)}}{ \prod\limits_{j=1}^{m-1}(1-\frac{1}{ q^{2dj}})} 
\end{eqnarray}
Again $m_{\phi}(g)=1$ and the  factor of $P_{|g|}(g)$ that corresponds to $\phi$ is 
$\frac{q^{-2d}}{ (1-\frac{1}{ q^{2d}})}.$  Likewise  $m_{\phi}(h)=m-1,$ and the factor of
$P_{|h|}(\lbrace h\rbrace )$ that corresponds to $\phi$ is $\frac{q^{-2d(m-1)}}{ \prod\limits_{j=1}^{m-1}(1-\frac{1}{ q^{2dj}})} $. Again these two expressions are precisely factors on the right side of
(\ref{AB2}).
\par
Finally, if $\phi$ does not divide $g$, then $m_{\phi}(g)=0$ and $m_{\phi}(f)=m_{\phi}(h).$
In this case,  the  factor of $P_{|g|}(g)$ that corresponds to $\phi$ is $1$,
and the factor of $P_{|h|}(h)$ that corresponds to $\phi$ is  exactly the same as the factor 
$P_{n}(\lbrace f\rbrace )$ that corresponds to $\phi$.

\end{proof}

\section{Estimating $E_{n}({\bf X}).$}

We know have all the tools necessary to prove the main result:
\begin{theorem}
\label{main}
 $\log E_{n}({\bf V})=n\log q-\log n+o_{q}(\log n).$
\end{theorem}
\begin{proof}
By Corollary  \ref{finalreduction}, it suffices to prove that
 $\log E_{n}({\bf X})=n\log q-\log n+o_{q}(\log n).$
 Recall the factorizations $f=\pi(f)h$, and define ${\cal G}_{n}=\lbrace g: g=\pi(f)$ for some $f\in \Omega_{n}\rbrace.$  Then 
\begin{align}E_{n}({\bf  X})=\sum\limits_{f\in \Omega_{n}}{\bf X}(\lbrace f\rbrace)P_{n}(\lbrace f\rbrace)&\\
=\sum\limits_{g\in {\cal G}_{n}}{\bf X}(\lbrace g\rbrace)
\sum\limits_{\lbrace h: \pi(gh)=g\rbrace}P_{n}(gh).
\end{align}
By Lemma \ref{sub}, this is less than or equal to
\begin{equation}
\label{sum1}
\sum\limits_{g\in {\cal G}_{n}}{\bf X}(\lbrace g\rbrace)
P_{|g|}(\lbrace g\rbrace)\sum\limits_{\lbrace h: h=f/g\text{ for some }g\in {\cal G}_{n}\rbrace}P_{|f/g|}(\lbrace h\rbrace).
\end{equation}
The inner sum is bounded by 1 since $P_{|f/g|}$ is a probability measure. 
Hence 
\begin{equation}
\label{dehed}
E_{n}({\bf X})\leq \sum\limits_{g\in {\cal G}_{n}}{\bf X}(\lbrace g\rbrace)
P_{|g|}(\lbrace g\rbrace).
\end{equation}
To estimate the sum in(\ref{dehed}), we need
 an upper bound for $P_{|g|}(\lbrace g\rbrace)$. 
Note that $|U(1,q^{d})|=q^{d}+1$ and $|GL(1,q^{2d})|=q^{2d}-1$ for all $d$.
Recall that, for $g\in {\cal G}_{n}$, we have $m_{\phi}(g)\leq 1$ for all $\phi\in {\cal  I}.$
Therefore, 
by Theorem \ref{charpoly}, we have
\begin{align}
&P_{|g|}(\lbrace g\rbrace)=\\
&q^{-|g|}\prod\limits_{\lbrace \phi\in {\cal J} :m_{\phi}(g)=1\rbrace}
\frac{1}{ 1+\frac{1}{ q^{|\phi|}}}
\prod\limits_{\lbrace \theta\in {\cal K}_{+} :m_{\theta}(g)=1\rbrace}
\frac{1}{ 1-\frac{1}{ q^{2|\phi|}}}\\
&\leq q^{-|g|}\prod\limits_{d=1}^{\infty}\frac{1}{ {1-\frac{1}{ q^{2d} } }}
\leq 2q^{-|g|}.
\end{align}

Thus 
\begin{equation}
E_{n}({\bf X})\leq 2\sum\limits_{g\in {\cal G}_{n}}q^{-|g|}{\bf X}(\lbrace g\rbrace)=
2\sum\limits_{m=1}^{n}q^{-m}\sum\limits_{\lbrace g\in{\cal G}_{n}: |g|=m\rbrace}{\bf X}(g).
\end{equation}

Factor each $g\in{\cal G}_{n}$ as $g=g_1g_2,$ where $g_1$ and $g_2$ respectively 
are the products of the irreducible factors in ${\cal J}$ and ${\cal K}$:
\begin{eqnarray}
g_1=\prod\limits_{\lbrace \phi\in {\cal J} :m_{\phi}(g)=1\rbrace}\phi&\\
g_2=\prod\limits_{\lbrace \theta\in {\cal K}_{+} :m_{\theta}(g)=1\rbrace}\theta\tilde{\theta}.
\end{eqnarray}
We certainly have 
\begin{equation}
{\bf X}(g)=LCM({\bf X}_{1}(g_{1}),{\bf X}_{2}(g_{2}))\leq {\bf X}_{1}(g_{1}){\bf X}_{2}(g_{2}),
\end{equation}
so
\begin{equation}
\label{crud}
E_{n}({\bf X})\leq 
2\sum\limits_{m=1}^{n}q^{-m}\sum\limits_{\lbrace g\in{\cal G}_{n}: |g|=m\rbrace }
{\bf X}_{1}(g_{1}){\bf X}_{2}(g_{2}).
\end{equation}
The degrees of of the irreducible factors of $g_{1}$ form a partition of the integer $|g_{1}|$ into distinct odd parts.  Let ${\cal S}_{1}(\pi)$ be the set of $g_{1}$'s with partition $\pi$.
Similarly,the degrees of the factors of $g_{2}$ from $\theta\in {\cal K}_{+}$ form a partition of $s:=\frac{|g_{2}|}{ 2}$ into distinct parts, and we let $S_{2}(\lambda)$ be the set of $g_{2}$'s with partition $\lambda.$ Using the notation $Q_{s}$ for the set of all partitions of $s$ into distinct parts, and
$O_{b}$ for the set of all partitions of $b$ into  distinct odd parts, we get

\begin{eqnarray}
\label{piece1}
\sum\limits_{\lbrace g\in{\cal G}_{n}: |g|=m\rbrace}{\bf X}_{1}(g_{1}){\bf X}_{2}(g_{2})=
\sum\limits_{s=1}^{\lfloor m/2\rfloor}\sum\limits_{\pi\in O_{m-2s}}
\sum\limits_{\lambda\in Q_{s}}|S_{1}(\pi)||S_{2}(\lambda)|{\bf X}_{1}(\pi){\bf X}_{2}(\lambda).
\end{eqnarray}

 Using the inequalites 
(\ref{Kd}) and (\ref{Jd}), we get
\begin{equation}
\label{piece2}
|S_{1}(\pi)|\leq \frac{q^{|\pi|}}{\pi_{1}\pi_{2}\cdots}
\end{equation}
and 
\begin{equation}
\label{piece3}
|S_{2}(\lambda)|\leq \frac{q^{2|\lambda|}}{\lambda_{1}\lambda_{2}\cdots}
\end{equation}
(where $\pi_{1},\pi_{2},\dots $ are the parts of $\pi$ and similarly for$\lambda$).
Putting (\ref{piece1}), (\ref{piece2}), and  (\ref{piece3}) back into the right side of
(\ref{crud}), we get

\begin{eqnarray}
\label{precesaro}
E_{n}({\bf X})\leq 
2\sum\limits_{m=1}^{n}
\left(\sum\limits_{s=1}^{\lfloor m/2\rfloor}\sum\limits_{\lambda\in O_{m-2s}}
\frac{{\bf X}_{1}(\pi)}{ \pi_{1}\pi_{2}\dots }
\sum\limits_{\lambda\in Q_{s}}\frac{{\bf X}_{2}(\lambda)}{ \lambda_{1}\lambda_{2}\dots }\right).
\end{eqnarray}
Let $\sigma_{2}(s)=\sum\limits_{\lambda\in Q_{s}}\frac{{\bf X}_{2}(\lambda)}{ \lambda_{1}\lambda_{2}\dots }$
be the innermost sum.
This sum was estimated  by Stong at the end of \cite{St1}. The conclusion was that, as $s\rightarrow\infty$,
\begin{equation}
\label{sig1bd}
\sigma_{2}(s)\leq \frac{(q^{2})^{s+o(\log s)}}{ s}.
\end{equation}
For any positve integer $b$ define
$\sigma_{1}(b)=\sum\limits_{\pi\in O_{b}}
\frac{{\bf X}_{1}(\pi)}{ \pi_{1}\pi_{2}\cdots }.$
We show next that it is sufficient to prove that 
\begin{equation}
\label{sigb}
\sigma_{1}(b)\leq \frac{q^{b+o(\log b)}}{ b}.
\end{equation}
Assume for now that (\ref{sigb}) holds. (It will be verified afterwards.)
For integers $k$ let $(k)^{+}=\max(k,1).$
Use the partial fraction decomposition $
\frac{1}{ s(m-2s)}=\frac{1}{ms}+\frac{2}{m(m-2s)}
$ so that
\begin{equation}
\sum\limits_{s=1}^{\lfloor m/2\rfloor}\frac{1}{s}\frac{1}{(m-2s)^{+}} =O(\frac{\log m}{m}).
\end{equation}
Then inside the parentheses of (\ref{precesaro}) we have
\begin{align}
&\sum\limits_{s=1}^{\lfloor m/2\rfloor }\sigma_{1}(m-2s)\sigma_{2}(s)=\\
&\sum\limits_{s=1}^{\lfloor m/2\rfloor }\frac{q^{m-2s+o(\log(m-2s))}}{(m-2s)^{+}}\frac{q^{2s+o(\log s)}}{s}\\
=&q^{m+o(\log m)}\sum\limits_{s=1}^{\lfloor m/2\rfloor}\frac{1}{s}\frac{1}{(m-2s)^{+}}=\frac{q^{m+o(\log m)}}{m}
\end{align}
Note that $q^{m}/m$ is an increasing function of $m$. So  if we let $\omega=\lfloor \log n\rfloor,$then
we can easily finish estimating (\ref{precesaro}):
\begin{align}
\sum\limits_{m=1}^{n}\frac{q^{m}}{m}= 
\sum\limits_{m=1}^{n-\omega}\frac{q^{m}}{m} +\sum\limits_{m=n-\omega +1}^{n}\frac{q^{m}}{m} \\
\leq (n-\omega)\frac{q^{n-\omega}}{n-\omega}+\omega \frac{q^{n}}{n}\\
=\frac{q^{n+o(\log n)}}{n}.
\end{align}

\par
To complete the proof of Theorem \ref{main}, all that remains is to prove that $\sigma_{1}(b)=
\frac{q^{b+o(\log b)}}{b}.$
The sum $\sigma_{1}$ is somewhat similar to $\sigma_{2}$, and we'll see that it can be estimated by
 techniques similar to those that Stong used in estimating $\sigma_{2}.$
The cyclotomic polynomials satisfy a simple identity:if $\pi_{i}$ is odd, then
\begin{equation}
q^{\pi_{i}}+1=\frac{q^{2\pi_{i}}-1}{ q^{\pi_{i}}-1}=\prod\limits_{d |\pi_{i}}\Phi_{2d}(q).
\end{equation}
Define
\begin{itemize}
\item
 ${ \Lambda}={\Lambda }(\pi)=\bigl\lbrace d:$ for some $i$, $d$ divides $\pi_{i}\bigr\rbrace$.
\item
$\nu_{d}(\pi)=\sum\limits_{k\equiv 0 (d)}c_{k}(\pi)=$ the number of parts that are multiples of $d$, and
\item$
w_{d}(\pi)=\max(0,\nu_{d}-1).
$
\end{itemize}
 Then 
\begin{align}
\label{prodqs}
LCM(q^{\pi_{1}}+1,q^{\pi_{2}}+1,\dots ) \leq
\prod\limits_{d\in {\Lambda}}\Phi_{2d}&\\
\label{numerator}
= \frac{\prod\limits_{i}(q^{\pi_{i}}+1)}{ 
\prod\limits_{d}\Phi_{2d}^{w_{d}}}.
\end{align}
If $\pi$ is a partion of $b$ into  distinct parts, then for the numerator of (\ref{numerator}) we have 
\begin{equation}
\prod\limits_{i}(q^{\pi_{i}}+1)=q^{b}\prod\limits_{i}(1-\frac{1}{ q^{\pi_{i}}})< q^{b}\prod\limits_{i=1}^{\infty}(1+\frac{1}{ q^{i}})< 4q^{b}.
\end{equation}
An upper bound is obtained if, in the denominator of (\ref{numerator}), we restrict $d$  to a finite set of primes.
For any $i$, let  $p_{i}$ denote the $i$'th prime; $p_{1}=2,p_{2}=3,\dots .$
Given a positive integer $\xi$, let ${\cal P}={\cal P}(\xi)=
\lbrace p_{i}: \xi\leq i\leq  e^{\xi}\rbrace=\lbrace p_{\xi},p_{\xi+1},\dots ,p_{\lfloor e^{\xi}\rfloor}\rbrace.$  Let $\kappa_{\xi}=\prod\limits_{p\in {\cal P}}\Phi_{2p}(q).$ Then, for any $\pi\in O_{b},$

\begin{equation}
\label{prodqs2}
LCM(q^{\pi_{1}}+1,q^{\pi_{2}}+1,\dots ) 
\leq \frac{4\kappa_{\xi}q^{b}}{
\prod\limits_{p\in {\cal P}}\Phi_{2p}^{\nu_{p}}}.
\end{equation}

Define
$$G(k)=\begin{cases}
 \prod\limits_{\lbrace p: p\in {\cal P} \text{ and } p|k\rbrace}\frac{1}{ \Phi_{2p}(q)},&
 \text{if $k$ is divisible by at least one prime in ${\cal P}$}\\
1&\text{else,}
\end{cases}$$

For any partition $\pi$, let 
$z_{\pi}=\frac{1}{ \prod\limits_{i=1}^{\infty}c_{i}!i^{c_{i}}}$,
where $c_{i}=c_{i}(\pi )$ 
is the number of parts of size $i$ that $\pi$ has. 
Thus $z_{\pi}=\frac{1}{ \pi_{1}\pi_{2}\cdots }$ for $\pi\in O_{b}$.
We get an upper bound for $\sigma_{1}(b)$ if we sum over {\sl all } partitions of $b$ 
(not just those in $O_{b}$). Hence
 and
from (\ref{prodqs}) we have
\begin{align}
\label{sig1bound}
\sigma_{1}(b)\leq  &4\kappa_{\xi}q^{b}\sum\limits_{\pi\vdash b}
\frac{z_{\pi}}{\prod\limits_{p\in {\cal P}}\Phi_{2p}^{\nu_{p}}}\\
&=4\kappa_{\xi}q^{b}\sum\limits_{\pi\vdash b}z_{\pi}\prod\limits_{k=1}^{\infty}G(k)^{c_{k}}
\end{align}

In the well-known cycle index identity 
\begin{equation}
1+\sum\limits_{b=1}^{\infty}
\sum\limits_{\pi\vdash b}z_{\pi}\prod\limits_{k}x_{k}^{c_{k}} z^{b}\
=\exp\left(\sum\limits_{k=1}^{\infty}\frac{x_{k}z^{k}}{ k}\right),
\end{equation}
we can make the substitutions $x_{k}=G(k), k=1,2,\dots$  to get
\begin{equation}
\label{upper2}
\sigma_{1}(b)\leq 4\kappa_{\xi}q^{b}
\lbrack\hskip-1pt\lbrack z^{b}\rbrack\hskip-1pt\rbrack 
\exp\left(\sum\limits_{k=1}^{\infty}\frac{G(k)}{ k}z^{k}\right)
\end{equation}

Following Stong,  we
note that the funtion $G(k)$ is a periodic function of $k$ with period 
$N=\prod\limits_{p\in {\cal P}}p$. Hence we have the Fourier expansion
\begin{equation}
G(k)=a_{0}+\sum\limits_{\ell=1}^{N-1}a_{\ell}\omega^{\ell k},
\end{equation}

where
$\omega=e^{2\pi i/N}$ and the $a_{\ell}$'s are the Fourier coefficients:
\begin{equation}
\label{Fourier}
a_{\ell}=\frac{1}{ N}\sum\limits_{v=0}^{N-1}G(v)\omega^{-\ell v}.
\end{equation}

Thus 
\begin{eqnarray}
\label{sing}
\exp\left(\sum\limits_{k=1}^{\infty}\frac{G(k)}{ k}z^{k}\right)
=
\exp\left(\sum\limits_{\ell=0}^{N-1}a_{\ell}
\sum\limits_{k=1}^{\infty}\frac{(\omega^{\ell}z)^{k}}{ k}\right)&\\
\label{sing2}
=(1-z)^{-a_{0}}\prod\limits_{\ell=1}^{N-1}(1-\omega^{\ell}z)^{-a_{\ell}}
\end{eqnarray}
Let $\alpha=\prod\limits_{j=1}^{N-1}(1-\omega^{j})^{-a_{j}}.$
Because $G(k)\geq 0$ for all $k$, it is clear from (\ref{Fourier}) that
$|a_{0}|>|a_{j}|  $ for all $j>0$. Hence the coefficient of $z^{n}$ in
(\ref{sing}) and (\ref{sing2}) is asymptotic to
\begin{equation}
\alpha \lbrack\hskip-1pt\lbrack z^{n}\rbrack\hskip-1pt\rbrack (1-z)^{-a_{0}}=
O( \frac{1}{n^{1-a_{0}}})
\end{equation}
It therefore suffices to verify that $a_{0}$ can be made arbitrarily small by choosing $\xi$ sufficiently large. 

\par
Note that, for odd primes $p$,  $\Phi_{2p}=\frac{q+1}{q^{p}+1}.$ 
 Hence
\begin{eqnarray}
\label{fbound}
G(k)\leq \begin{cases} 1& \text{if } $g.c.d.(k,N)=1$\\
    \frac{q+1}{ q^{p_{\xi}}+1}
&\text{else}
\end{cases}
\end{eqnarray}
Given $\epsilon>0$, choose $\xi $ large enough so that we also have
\begin{equation}\frac{q+1}{ q^{p_{\xi}}+1}<\epsilon/2.\end{equation}
\par

Let $R_{\xi}=\lbrace k: gcd(k,N)=1 \text{ and } k\leq N\rbrace$
By inclusion-exclusion, 
$|R_{\xi}|=\prod\limits_{i=\xi}^{\lfloor e^{\xi}\rfloor}(1-\frac{1}{ p_{i}})N$.
By the prime number theorem $p_{i}\sim i\log i$  and consequently
$\prod\limits_{i=\xi}^{\lfloor e^{\xi}\rfloor}(1-\frac{1}{ p_{i}})=o(1) $ as $\xi\rightarrow\infty.$  
 We can therefore also choose
$\xi$ large enough so that 
  $|R_{\xi}|\leq \frac{\epsilon}{ 2}N.$ 
But then
\begin{eqnarray}
a_{0}=\frac{1}{ N}\sum\limits_{k=1}^{N-1}G(k)\leq \frac{|R_{\xi}|}{N}+\frac{p_{\xi}+1}{ q^{p_{\xi}}+1}<\epsilon.
\end{eqnarray}
\end{proof}


\begin{thebibliography}{9}
\bibitem{Berl} Elwyn R. Berlekamp, Algebraic Coding Theory,2nd Ed., Aegean Park Press,  ISBN: 0-89412-063-8
\bibitem{Fulman1}
Jason Fulman, Cycle indices for the finite classical groups. {\sl J. Group Theory} {\bf 2} (1999), no. 3, 251--289, arXiv:math.GR/9712239
\bibitem{FPN} 
Jason Fulman, Peter M. Neumann, and Cheryl E. Praeger,  
A generating function approach to the enumeration of matrices in classical groups over finite fields.
{\sl Mem. Amer. Math. Soc. } {\bf 176 } (2005), no. 830.
\bibitem{Ful3} Jason Fulman, Random matrix theory over finite fields. {\sl Bull. Amer. Math. Soc. (N.S.) } {\bf 39} (2002), no. 1, 51--85, arXiv:math.GR/0003195
\bibitem{Grove} Larry C.Grove, Classical groups and geometric algebra,
{\sl Graduate Studeies in Mathematics } {\bf 39}, American Mathematical Society (2002), ISBN 0-8218-2019-2.
\bibitem{Kung}Joseph P.S. Kung,  
The cycle structure of a linear transformation over a finite field. 
Linear Algebra Appl. 36 (1981), 141--155.	
\bibitem{Lidl} R. Lidl and H.Niederreiter, Introduction to finite fields  and their applications, Cambridge University Press (1994).
\bibitem{MeinGotz}Helmut Meyn and Werner  G\"otz,
Self-reciprocal polynomials over finite fields. 
{\sl SŽminaire Lotharingien de Combinatoire (Oberfranken, 1990), }82--90, 
{\tt Publ. Inst. Rech. Math. Av.,} {\bf 413}, 
Univ. Louis Pasteur, Strasbourg, 1990. 
\bibitem{St1}Richard Stong, The average order of a matrix. {\sl  J. Combin. Theory Ser. A }{ \bf 64} (1993), no. 2, 337--343.

\bibitem{St2} Richard Stong, 
Some asymptotic results on finite vector spaces. 
{\sl Adv. in Appl. Math.} {\bf 9} (1988), no. 2, 167--199.
\bibitem{YucasMullen} Joseph L. Yucas and Gary L. Mullen, Self-Reciprocal Irreducible Polynomials Over Finite Fields, 
{\sl Designs, Codes, and Cryptography} {\bf 33} (2004) 275--281.
\bibitem{Wall}G.E. Wall,  
On the conjugacy classes in the unitary, symplectic and orthogonal groups. 
{\sl J. Austral. Math. Soc.} {\bf 3} ( 1963) 1--62.
\end{thebibliography}
\end{document}